# THE OUTSIDE OF THE TEICHMÜLLER SPACE OF PUNCTURED TORI IN MASKIT'S EMBEDDING

JOUNI PARKKONEN

University of Jyväskylä

September 15, 1997

ABSTRACT. We show that for each cusp on the boundary of Maskit's embedding $\mathscr{M} \subset \mathbb{H}$ of the Teichmüller space of punctured tori there is a sequence of parameters in the complement of $\overline{\mathscr{M}}$ converging to the cusp such that the parameters correspond to discrete groups with elliptic elements. Using Tukia's version of Marden's isomorphism theorem we identify them as cusps on the boundary of certain deformation spaces of Koebe groups with a non-simply connected invariant component.

Let $G$ be a Kleinian group with a simply connected *invariant component* $\Omega_0$ of the *set of discontinuity* $\Omega(G)$ of $G$. $G$ is a *terminal b-group* of type $(p,n)$ if $\Omega_0/G$ is a Riemann surface of finite analytic type $(p,n)$ and if $(\Omega \setminus \Omega_0)/G$ is a collection of $2p - 2 + n$ thrice punctured spheres. We use the matrix notation for Möbius transformations, identifying the transformation

$$z \mapsto \frac{az+b}{cz+d}$$

with the matrix

$$\begin{pmatrix} a & b \\ c & d \end{pmatrix} \in \mathrm{PSL}(2, \mathbb{C}).$$

The *Maskit embedding* [14], [10] of the Teichmüller space $\mathscr{T}(p,n)$ of compact Riemann surfaces of genus $p$ with $n$ punctures is defined by first identifying $\mathscr{T}(p,n)$ with the *(quasiconformal) deformation space* $\mathscr{T}(G)$ of a terminal b-group,

$$\mathscr{T}(G) = \left\{ w \colon \widehat{\mathbb{C}} \to \widehat{\mathbb{C}} \text{ quasiconformal} : w \circ g \circ w^{-1} \in \mathrm{PSL}(2,\mathbb{C}) \text{ for all } g \in G \right\} / \sim,$$

where $w_1 \sim w_2$ if there is a Möbius transformation $A \in \mathrm{PSL}(2,\mathbb{C})$ so that

$$w_1 \circ g \circ w_1^{-1} = A \circ w_2 \circ g \circ w_2^{-1} \circ A^{-1} \quad \forall g \in G.$$

This deformation space is then, using traces and/or fixed points of a collection of group elements, analytically embedded in $\mathbb{C}^N$ where $N = \dim \mathscr{T}(G) = \dim \mathscr{T}(p,n)$.

Research supported by the Academy of Finland (Grant 34082) and Magnus Ehrnroothin säätiö of the Finnish Society of Sciences and Letters

Typeset by $\mathcal{A}_{\mathcal{M}}\mathcal{S}$-TEX

In this paper we study the one-dimensional case of once-punctured tori. For once-punctured tori the Maskit embedding rather simple to describe (See Kra [10], Section 6 for more details):

$\mathscr{T}(1,1)$ can be represented as the set $\mathscr{M}$ of parameters $\mu \in \mathbb{H}$ for which the group $G[\mu]$ generated by two transformations

$$S = \begin{pmatrix} 1 & 2 \\ 0 & 1 \end{pmatrix}$$

and

$$T[\mu] = \begin{pmatrix} -i\mu & -i \\ -i & 0 \end{pmatrix}$$

is a terminal b-group. It is useful to think of $G$ as the HNN extension (See Maskit [16], Section VII.D) of a torsion free triangle group $\langle S, \widetilde{S} \rangle$ by the loxodromic transformation $T[\mu]$, where

$$\widetilde{S} = T[\mu]^{-1} S^{-1} T[\mu] = \begin{pmatrix} 1 & 0 \\ -2 & 1 \end{pmatrix}.$$

The group $\langle S, \widetilde{S} \rangle$ is a Fuchsian group of the first kind that keeps $\mathbb{H}$ and $\mathbb{H}^*$ fixed. The quotient $\Omega(G[\mu])/G[\mu]$ consists of the disjoint union of a punctured torus $\Omega(G[\mu])/G[\mu]$, where $\Omega(G[\mu]) = \Omega_0[\mu]$ is the invariant component of $G[\mu]$, and a thrice punctured sphere. The punctures of the sphere correspond to the three parabolic conjugacy classes of $S$, $\widetilde{S}$ and

$$K = S\widetilde{S} = \begin{pmatrix} -3 & 2 \\ -2 & 1 \end{pmatrix}.$$

$K$ also corresponds to the puncture on the torus component of the quotient Riemann surface.

The set $\mathscr{M}$ and its boundary have been studied in detail by Keen and Series in [7] and by Wright in the unpublished manuscript [25]. In this note we consider discrete groups generated by $S$ and $T[\mu]$ for parameters $\mu \in \mathbb{C} \setminus \overline{\mathscr{M}}$. We can restrict our attention to the case $\mu \in \overline{\mathbb{H}} \setminus \mathscr{M}$, as a simple calculation shows that for all $\mu \in \mathbb{C}$

$$G[\mu] = E G[-\mu] E,$$

for $E(z) = -z$, $E^2 = 1$. We consider the following questions:

(1) For which parameters $\mu \in \mathbb{C} \setminus \overline{\mathscr{M}}$ is the group $G[\mu]$ discrete?
(2) If $G[\mu]$ is discrete, what is the geometry of the quotients $\Omega(G[\mu])/G[\mu]$ and $\mathbb{H}^3/G[\mu]$?

Using the circle chain methods of Keen and Series [7] and Wright [25] and Maskit's second combination theorem [17] we prove a local result: (see Section 1 for the definitions)

**6.5. Theorem.** *On the extension $\mathscr{P}_{p/q}^+$ of each rational pleating ray there is an open neighborhood $U$ of the cusp $\mu_{p/q}$ on $\mathscr{P}_{p/q}^+$ such that the following holds: If $\mu \in \mathscr{P}_{p/q}^+ \cap U$ and $\operatorname{tr} W_{p/q} = 2\cos(\pi/n)$ for some $n \in \mathbb{N}$, $n \geq 2$, then*

$$G[\mu] = F *_{W_{r/s}}$$

*is a Kleinian group representing a thrice punctured sphere and a sphere with a puncture and two branch points of order $n$ on its ordinary set. If $\mu \in \mathscr{P}_{p/q}^+ \cap U$ is not of this form, $G[\mu]$ is not discrete.*

A similar observation for the Bers embedding of Teichmüller space is made by Kra and Maskit in [11], and for the Riley slice by Riley in the introduction to [21].

We do not have an estimate for the size of the neighborhood $U$, that is, we do not know how big $n$ has to be for a fixed rational $p/q$ in order that the group is discrete. For integral pleating rays we get all $n \geq 2$, and computer experiments suggest that this should be true for the extension of any rational pleating ray.

In Section 6 we also show (Theorem 6.6) that the discrete groups with elliptic elements that we find in the complement of $\overline{\mathscr{M}}$ are boundary groups of deformation spaces $\mathscr{M}_n$ of certain Kleinian groups representing a punctured torus on their invariant component.

**6.6. Theorem.** *Let $\tau_{p'/q}$ be the boundary point of $\mathscr{M}_n$ on the ray $\mathscr{P}_{p'/q}$. Then the group $G[\mu_{p/q}(n)]$ is conjugate to $G_n[\tau_{p'/q}]$ by a Möbius transformation.*

The quotient 3-orbifold $M = \mathbb{H}^3/G$ of these groups can be described topologically/geometrically as follows: Let $X$ be a punctured torus. M is obtained from $X \times [0,1]$ by adding a singular two-handle $\mathbb{D}_n \times [0,1]$, where $\mathbb{D}_n$ is a disk with metric cone singularity of order $n$ in the origin. Note that one can construct these groups for any $n \geq 2$, however, the method of proof of Theorem 6.5 does not allow us to conclude that the boundary groups are realized on the extended rays for small $n$. The deformation spaces $\mathscr{M}_n$ are treated in [1] and [20].

## Plan of the paper

In Section 1 we sketch the construction of rational pleating rays in $\mathscr{M}$ due to Keen and Series [6],[7], and define extended rays. In Section 2 we consider disk-preserving subgroups of $G[\mu]$ for parameters $\mu$ that are contained in an extended rational pleating ray, and derive a necessary condition for the discreteness of $G[\mu]$. Sections 3 and 4 are concerned with modifying the definitions of pleating rays and circle chains of Wright [25] to the case of groups with elliptic elements. In the end of Section 4 we prove Theorems 6.5 and 6.6 for $p/q \in \mathbb{Z}$. Sections 5 and 6 contain the proofs of the main results. The proof of Theorem 6.5 is split in to two parts, the first part (Theorem 5.1) is an application of Maskit's second combination theorem [17] that reduces the proof of Theorem 6.5 to finding a circle chain with a number of special properties. This chain is then constructed in Section 6.

## Acknowledgements

I am grateful to Linda Keen for suggesting the problem "What happens to pleating rays on the other side of the boundary of $\overline{\mathscr{M}}$?" at the special session in honor of Bernard Maskit's 60th birthday at the AMS meeting at Hartford, Connecticut in March 1995. The work presented in this paper is the result of my attempt to understand this question.

Part of the research for this paper was done during my stay at the Mathematics Department and IMS at Stony Brook in 1994-1995 and at Institut Henri Poincaré, Centre Emile Borel in 1996. I would like to thank these institutes and the organizers of the Semester on Low Dimensional Geometry and Topology (CEB/IHP) for their hospitality.

I would also like to thank Curtis T. McMullen and David J. Wright for the computer programs used to produce the limit sets of Kleinian groups and the picture of the deformation space. Figures 2 and 8 were computed using McMullen's program, and Figures 3 and 5 using Wright's program.

## NOTATION

For the convenience of the reader we list here some of the most frequently used notations and conventions used in this article:

| | |
|---|---|
| $\widehat{\mathbb{C}}$ | The extended complex plane $\mathbb{C} \cup \{\infty\}$, |
| $\mathbb{H}$ | the upper half plane $\{z \in \mathbb{C} : \operatorname{Im} z > 0\}$, |
| $\mathbb{H}^3$ | the upper half space $\mathbb{C} \times \{t \in \mathbb{R} : t > 0\}$, |
| $\operatorname{PSL}(2,\mathbb{C})$ | the projective special linear group of complex $2 \times 2$ matrices identified with the group of Möbius transformations of $\widehat{\mathbb{C}}$, |
| $\Lambda(G)$ | The limit set of the Kleinian group $G$, |
| $\Omega(G)$ | The set of discontinuity of the Kleinian group $G$. |

If $X_1$ and $X_2$ are Möbius transformations, we denote the group they generate by $\langle X_1, X_2 \rangle$. Similarly, $\langle G, X \rangle$ denotes the group generated by the elements of the group $G$ and the Möbius transformation $X$. $W_{p/q}[X_1, X_2]$ denotes a special word in the transformations $X_1$ and $X_2$, see Section 1 for details and references.

We also use the following notations: Defined:

| | | |
|---|---|---|
| $\mathscr{M}$ | Maskit's embedding of $\mathscr{T}(1,1)$ | Introd. |
| $\mathscr{M}_n$ | 'Maskit-type' embedding of the deformation space of terminal Koebe groups of type $(1,1)$ | Section 3 |
| $\mathscr{P}_{p/q}$ | The $p/q$ pleating ray in $\mathscr{M}$ | Section 1 |
| $\mathscr{P}_{p/q}^+$ | The extended $p/q$ pleating ray for $\mathscr{M}$ | Section 1 |
| $\mathscr{P}_{p/q}^n$ | The $p/q$ pleating ray in $\mathscr{M}_n$ | Section 4 |
| $\mu_{p/q}$ | The unique finite point in $\overline{\mathscr{P}_{p/q}} \cap \mathscr{M}$ | Section 1 |
| $\mu_{p/q}(n)$ | A point in $\mathscr{P}_{p/q}^+$ such that $\lvert \operatorname{tr} W_{p/q} \rvert = 2\cos(\pi/n)$ | Section 1 |
| $\tau_{p/q}$ | The unique finite point in $\mathscr{P}_{p/q} \cap \mathscr{M}_n$ | Section 1 |

*Further conventions.* $\mu$ denotes a parameter in the complex plane containing $\mathscr{M}$. $\tau$ denotes a parameter in the complex plane containing $\mathscr{M}_n$, where $n = 2, 3, 4, \ldots$. If a Möbius transformation $X$ or a group $G$ depends on a parameter $\xi$, this dependence is denoted by $X[\xi]$ or $G[\xi]$, respectively.

## 1. EXTENDING RATIONAL PLEATING RAYS OF $\mathscr{M}$

Keen and Series studied the Maskit embedding in [7]. They related the parameters in $\mathscr{M}$ to the geometry of the quotient of $\mathbb{H}^3$. In this section we sketch the definition of a rational pleating ray and define an extension of rational pleating rays into the complement of $\mathscr{M}$.

Free homotopy classes of simple closed curves on a punctured torus, and the corresponding group elements, special words in the generators $S$ and $T[\mu]$, can be enumerated by $\mathbb{Q} \cup \{\infty\}$. The special words

$$W_{p/q} = W_{p/q}[\mu] = W_{p/q}[S, T[\mu]]$$

can be defined inductively using the Farey sequence as follows (see [7], Osborne and Zieschang [19], or Cohen, Metzler and Zimmermann [3] for details): Two rational numbers $p/q$ and $p'/q'$ are called *neighbors* if

$$pq' - p'q = \pm 1.$$

First set
$$W_{1/0} = W_\infty = S^{-1} \quad \text{and} \quad W_{0/1} = W_0 = T.$$

If $a/b < c/d$ are neighbors, set

$$W_{(a+c)/(b+d)} = W_{c/d} W_{a/b}.$$

For each rational number $p/q$ we consider the set $\widetilde{\mathcal{V}}_{p/q}$ where the word $W_{p/q}[\mu]$ has real trace. We call any curve in the free homotopy class of the simple closed geodesic $\gamma_{p/q}$ on $\Omega_0[\mu]/G[\mu]$ determined by $W_{p/q}$ a $p/q$-*curve*.

Keen and Series [7] used the geometry of the boundary of the convex core of the 3-manifold $M = \mathbb{H}^3/G[\mu]$ to parameterize $\mathscr{M}$: Let $\partial \mathscr{C}_0 = \mathscr{C}_0[\mu]$ be the component of the hyperbolic convex hull of the limit set $\Lambda(G[\mu])$ that lies above the invariant component $\Omega_0$. The limit set, and thus its convex hull is $G[\mu]$ invariant, so $G[\mu]$ acts on $\partial \mathscr{C}_0$. The quotient surface $\partial \mathscr{C}_0/G[\mu]$ is a pleated surface in the sense of Thurston [23], topologically it is a punctured torus. It inherits a geometric structure from the embedding of $\partial \mathscr{C}_0$ in $\mathbb{H}^3$, in which it is planar outside a geodesic lamination.

The $p/q$ (rational) pleating ray $\mathscr{P}_{p/q}$ is defined in [7] as the locus of $\mu \in \mathscr{M}$ for which the surface $\partial \mathscr{C}_0/G[\mu]$ is pleated along the $p/q$ curve. The space $\mathscr{PML}$ of projective measured laminations on a punctured torus can be identified with $\widehat{\mathbb{R}} = \mathbb{R} \cup \{\infty\}$. Keen and Series [9] showed that the map pl: $\mathscr{T}(G) \to \widehat{\mathbb{R}}$ defined by mapping a parameter $\mu$ to the pleating lamination of $\partial \mathscr{C}_0/G[\mu]$ is continuous.

A rational pleating ray $\mathscr{P}_{p/q}$ can be characterized as the unique asymptotically vertical connected component of

$$\widetilde{\mathcal{V}}_{p/q} \cap \mathscr{M}.$$

The pleating ray $\mathscr{P}_{p/q}$ is homeomorphic to the real line and it ends at a *cusp* $\mu_{p/q}$ on $\partial \mathscr{M}$. In this note we consider an extension of the ray across the boundary of $\mathscr{M}$:

**1.1. Definition.** $\mathscr{P}^+_{p/q}$, the connected component of

$$\overline{\mathscr{P}_{p/q}} \cup \left\{ \mu \in \widetilde{\mathcal{V}}_{p/q} \cap (\mathbb{H} \setminus \overline{\mathscr{M}}) : |\operatorname{tr} W_{p/q}| > 0 \right\}$$

containing $\mathscr{P}_{p/q}$ is the *extended $p/q$ ray*.

*1.2. Remark.* The trace function may have a critical point close to the boundary even at $\mu_{p/q}$, that is, $\mathscr{P}^+_{p/q}$ might have a branching point. The somewhat awkward definition of the extended ray is used in this form so that we only get the correct ray corresponding to the pleating ray *inside $\mathscr{M}$*.

**1.3. Lemma.** *There are points $\mu_{p/q}(n) \in \mathscr{P}_{p/q}^+$ for all large values of $n$ such that $W_{p/q}[\mu_{p/q}(n)]$ is elliptic of order $n$, and*

$$\lim_{n \to \infty} \mu_{p/q}(n) = \mu_{p/q}.$$

*Proof.* The function $\mu \mapsto \operatorname{tr} W_{p/q}[\mu]$ is holomorphic. Thus, the image of a small neighborhood of $\mu_{p/q}$ covers a small neighborhood of $\operatorname{tr} W_{p/q} = \pm 2$, thus for large $n$ the values $\operatorname{tr} W_{p/q}[\mu] = \pm 2\cos(\pi/n)$ are attained. □

In most of what follows we aim at proving that that for large $n$ the groups $G[\mu_{p/q}(n)]$ are discrete. The method is based on a 'reconstruction' of $G[\mu]$ starting from a Fuchsian subgroup. We can only find a subgroup like this for rational rays, and do not know how to treat groups corresponding to parameters that are not contained in an extended rational ray. However, we have the following result as an application of Jørgensen's inequality [4]:

**1.4 Proposition.** *For each $\mu_{p/q}(n)$ there is an open neighborhood $U_{p/q}(n)$ such that if $\mu \in U_{p/q}(n)$, then $G[\mu]$ is not discrete.*

*Proof.* Consider the subgroup

$$H_{p/q}(n) = \langle K, W_{p/q}^n \rangle,$$

where

$$\begin{pmatrix} -3 & 2 \\ -2 & 1 \end{pmatrix}.$$

By Jørgensen's inequality [4] Lemma 1, we know that if

$$|\operatorname{tr}(KW_{p/q}[\mu]^n K^{-1} W_{p/q}[\mu]^{-n} - 2| < 1,$$

then either $H_{p/q}(n)$ is elementary or it is not discrete. If $mu = \mu_{p/q}(n)$, then $KW_{p/q}[\mu]^n K^{-1} W_{p/q}[\mu]^{-n}$ is the identity, and $H_{p/q}(n) = \langle K \rangle$ is discrete. If $W_{p/q}[\mu]$ is loxodromic, $H_{p/q}(n)$ cannot be elementary. Also, if $KW_{p/q}[\mu]^n K^{-1} W_{p/q}[\mu]^{-n}$ is elliptic, then it has to be of order 2, 3, 4, or 6 by the classification of elementary groups (see Maskit Section V.D). The combination of these observations gives the desired neighborhood. □

## 2. A necessary condition for discreteness of $G[\mu]$

In this section we restrict to a subgroup $F$ of $G[\mu]$ that preserves a circle, and find the values of $\mu \in \mathscr{P}_{p/q}^+$ for which this subgroup is discrete. In [7] Keen and Series consider the following subgroup of $G[\mu]$ for parameters $\mu \in \mathscr{M}$: Let $p/q$ and $r/s > p/q$ be Farey neighbors, and set

$$F = F_{p/q,\mu} = \langle W_{p/q}, W_{r/s}^{-1} W_{p/q}^{-1} W_{r/s} \rangle.$$

They show that if $\mu \in \mathscr{P}_{p/q}$, then $F_{p/q}$ is a Fuchsian group of the second kind, that represents a punctured cylinder with boundary geodesics of equal length. For parameters in the extension of the $p/q$ ray we have:

**2.1. Proposition.** *Let $\mu \in \mathscr{P}_{p/q}^+ \setminus \overline{\mathscr{M}}$. If $G[\mu]$ is discrete, then $|\operatorname{tr} W_{p/q}| = 2\cos(\pi/n)$ for some $n = 2, 3, \ldots$.*

The proof of the Proposition is divided into two Lemmas. The first one is an analog of Proposition A.1 of [7].

**2.2. Lemma.** *If $\mu \in \mathscr{P}_{p/q}^+$ and $\operatorname{tr} W_{p/q} = 2\cos(\pi/n)$ for some $n \in \mathbb{N}$, $n \geq 2$, then $F = F_{p/q,\mu}$ is a triangle group of signature $(n, n, \infty)$.*

*Proof.* The elements $W_{p/q}$, $W_{r/s}^{-1} W_{p/q}^{-1} W_{r/s}$ and their product $K$ all have real traces. If $n \geq 3$ it follows from Beardon [2] (Proof of Theorem 5.2.1) that $F$ preserves a disk $D$. The quadrilateral with vertices at the fixed points of $W_{p/q}$, $W_{r/s}^{-1} W_{p/q}^{-1} W_{r/s}$, $K$ and
$$\widetilde{K} = \widetilde{K}_{p/q} = W_{p/q} K W_{p/q}^{-1}$$
in $\overline{D}$ satisfy the conditions of Poincaré's theorem. If $n = 2$ we do not get a group that acts on a disc, but an infinite dihedral group acting on $\mathbb{C}$. $\square$

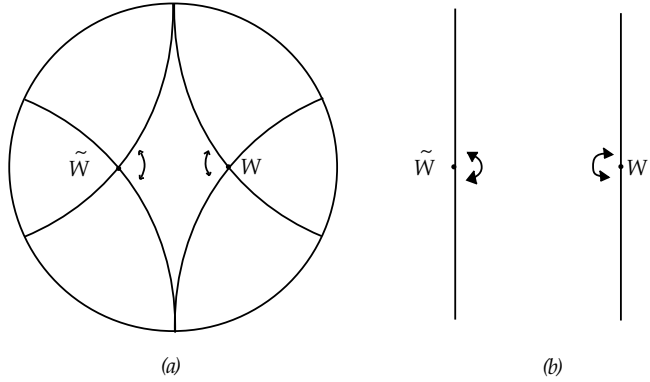

FIGURE 1. The Ford fundamental polygon of (a) a triangle group of signature $(n, n, \infty)$ acting in $\mathbb{D}$, and (b) the Euclidean triangle group of signature $(2, 2, \infty)$ acting in $\mathbb{C}$.

Let $G$ be a Kleinian group. An elliptic Möbius transformation $g \in G$ is a *primitive* element of $G$, if it is conjugate to a rotation by $2\pi/n$, $n \in \mathbb{Z}$, and if any elliptic element in $G$ with the same fixed points as $g$ is a power of $g$.

**2.3. Lemma.** *If $\mu \in \mathscr{P}_{p/q}^+$ and $\operatorname{tr} W_{p/q} \in (-2, 2) \setminus \{2\cos(\pi/n) : n \in \mathbb{N}, n \geq 2\}$, then $F = F_{p/q,\mu}$ is not discrete.*

*Proof.* If $W_{p/q}$ is an infinite order elliptic, the group is trivially non-discrete, so we can assume $W_{p/q}$ is a non-primitive elliptic. Let $W_{\text{prim}}$ be a primitive generator of $\langle W_{p/q} \rangle$. $F$ has parabolic elements, so any fundamental polygon has cusps extending to the circle at infinity. The commutator relation implies that the isometric circles of $W_{p/q}$ and $W_{r/s}^{-1} W_{p/q}^{-1} W_{r/s}$ are tangent at the fixed points of $K$ and $\widetilde{K}$. Clearly the isometric circles of $W_{\text{prim}}^{\pm 1}$ and $W_{r/s}^{-1} W_{\text{prim}}^{\mp 1} W_{r/s}$ intersect inside $D$. The Ford fundamental region of $F$ would thus be bounded away from the circle at infinity, which is impossible because there are parabolic elements in the group. $\square$

## 3. Koebe groups of type $(1, 1, \infty)$

In this Section we introduce parameter spaces similar to the Maskit embedding for a class of geometrically finite Koebe groups that represent a punctured torus on their invariant component. Koebe groups representing Riemann surfaces with cone points and punctures, and their deformation spaces were treated in [20]. We are able to simplify some of the expressions of [20] considerably by restricting to the case of the *punctured* torus.

These groups are introduced here for the following reason: In Theorems 4.8 and 6.5 we prove that the groups $G[\mu_{p/q}(n)]$, in the cases we can show that they are discrete, are HNN-extensions of triangle groups of signature $(n, n, \infty)$. In fact, in Theorems 4.11 and 6.6 we are able to show that $G[\mu_{p/q}(n)]$ is conjugate by a Möbius transformation to a group corresponding to a boundary point of one of the deformation spaces, $\mathscr{M}_n$ introduced in this Section.

The following definition is due to Maskit:

**3.1 Definition.** A Kleinian group $G$ that has an invariant component $\Omega_0 \subset \Omega(G)$ is called a *Koebe group* if any other component $\Delta' \subset \Omega(G) \setminus \Omega_0$ is a round disc. If the stabilizers of the disk components are hyperbolic triangle groups, the Koebe group is *terminal*.

*3.2. Remark.* This is a large class of groups, containing b-groups (simply connected invariant component) and (quasi-)Fuchsian groups (two invariant components). Here Koebe group usually means a Koebe group that is not a b-group.

The groups we consider are extensions of triangle groups: Let $F_n$, $2 \leq n \leq \infty$, be a triangle group of signature $(n, n, \infty)$. If $n > 2$ $F_n$ is a hyperbolic triangle group, $F_2$ is the infinite dihedral group acting on $\mathbb{C}$. Let $A$ and $B$ be *canonical generators* for $F_n$, that is, for $2 \leq n < \infty$, $A$ and $B$ are primitive elliptic elements in $F_n$ with the presentation

$$F_n = \langle A_n, B_n : A_n^n = B_n^n = \mathrm{id}, \ K_n = A_n B_n \text{ parabolic} \rangle,$$

and for $n = \infty$, $A_\infty$ and $B_\infty$ are primitive parabolics in $F_\infty$ with the presentation

$$F_\infty = \langle A_\infty, B_\infty : A_\infty, B_\infty, \text{ and } K_\infty = A_\infty B_\infty \text{ parabolic} \rangle.$$

We wish to construct a group $\langle F_n, C \rangle$ that realizes the following gluing construction: Let $X$ be a component of the quotient $\Omega(F_n)/F_n$. We cut out disks around the special points/punctures on $X$ corresponding to the generators $A_n$ and $B_n$, and glue the boundaries of the remaining surface together, thus producing a punctured torus. It is relatively easy that this is achieved by adding a new generator $C_n[\tau]$n satisfying

$$C_n[\tau]^{-1} A_n C_n[\tau] = B_n^{-1}.$$

There is a one complex parameter family of transformations $C_n[\tau]$ satisfying this condition. Set

$$G_n[\tau] = \langle F_n, C_n[\tau] \rangle.$$

If we choose a normalization for our groups we can be more explicit:

The extensions of the torsion-free triangle group $F_\infty$ are groups of the form

$$G[\mu] = G_\infty[\mu] = \langle S, \widetilde{S} \rangle *_{T[\mu]} = \langle S, T[\mu] \rangle.$$

It was shown in [10] that for $\operatorname{Im}\mu > 2$, $G[\mu]$ is a terminal b-group uniformizing a punctured torus. The Maskit embedding $\mathscr{M}$ is the open set of parameters $\mu \in \mathbb{H}$, containing $\{\mu|\operatorname{Im}\mu > 2\}$ for which $G[\mu]$ is a terminal b-group. For these parameters

$$G[\mu] = F_\infty * C_\infty[\mu].$$

Let $3 \leq n < \infty$. After normalization we can assume that $F_n$ is generated by two elliptic transformations

$$A_n = \begin{pmatrix} e^{-i\pi/n} & 0 \\ 0 & e^{i\pi/n} \end{pmatrix},$$

and

$$B_n = \begin{pmatrix} i\sin(\pi/n)\cosh d_n - \cos(\pi/n) & -i\sin(\pi/n)\sinh d_n \\ i\sin(\pi/n)\sinh d_n & -i\sin(\pi/n)\cosh d_n - \cos(\pi/n) \end{pmatrix}$$

$$= \begin{pmatrix} 2i/\sin(\pi/n) - e^{i\pi/n} & -2i\cot(\pi/n) \\ 2i\cot(\pi/n) & -2i/sin(\pi/n) - e^{-i\pi/n} \end{pmatrix}.$$

where

$$d_n = \operatorname{arcosh}\frac{\cos^2 \pi/n + 1}{\sin^2 \pi/n} = 2\operatorname{arcosh}\frac{1}{\sin(\pi/n)}$$

is the hyperbolic distance between the fixed points of $A$ and $B$ in $\mathbb{D}$. Now

$$C_n[\tau] = \begin{pmatrix} \tau\sinh(d_n/2) & -\tau\cosh(d_n/2) \\ \frac{1}{\tau}\cosh(d_n/2) & -\frac{1}{\tau}\sinh(d_n/2) \end{pmatrix}$$

$$= \begin{pmatrix} \tau\cot(\pi/n) & -\tau/\sin(\pi/n) \\ 1/(\tau\sin(\pi/n)) & -\cot(\pi/n)/\tau \end{pmatrix}.$$

It was shown in [20] that for

$$|\tau| > \coth(d_n/4) = \frac{1 + \sin(\pi/n)}{\cos(\pi/n)},$$

the group

$$G_n[\tau] = \langle F_n, C_n[\tau] \rangle$$

is a Koebe group that represents a punctured torus on its invariant component, and that the map $G_n[\tau] \mapsto \tau^2$ is a global coordinate for the deformation space of any fixed Koebe group $G_n[\tau_0]$.

**3.3. Definition.** We denote the image of the deformation space of $G_n[\tau_0]$ by $\mathscr{M}_n$.

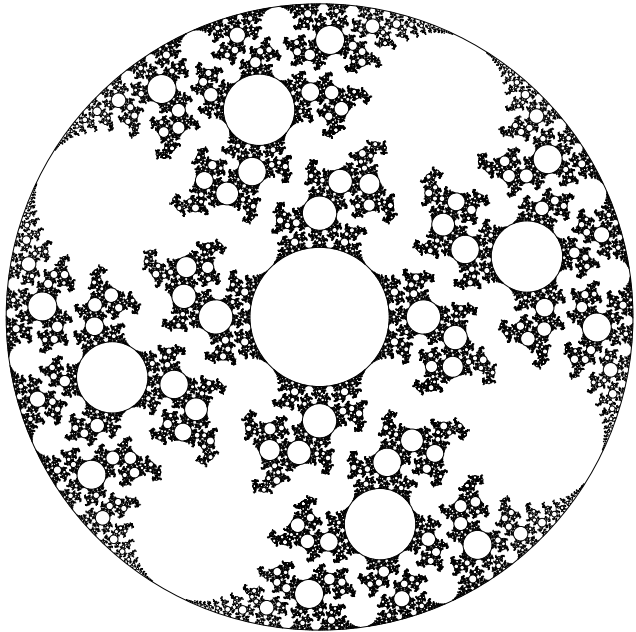

FIGURE 2. The limit set of a Kleinian group from $\mathcal{M}_4$ corresponding to the parameter $\tau = 2.1 + 0.3\ i$.

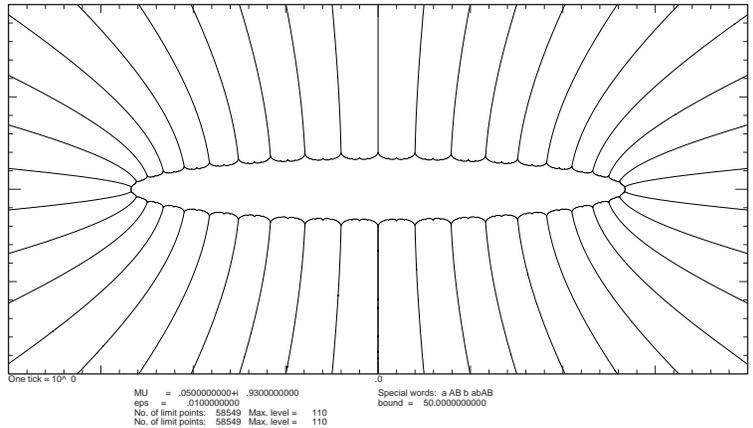

FIGURE 3. This picture shows a 2 to 1 covering of $\mathcal{M}_{21}$ with the lifts of integral pleating rays $\mathcal{P}_m^{21}$, $m = 0, 1, \ldots, 20$. See Remark 3.4(2) for more details.

The extensions of the infinite dihedral group $F_2$ are not that important for the rest of the paper, but it is interesting to note the same constructions as in the case of hyperbolic triangle groups can be performed. $F_2$ can be generated by two elliptics of order 2 (half turns)

$$A_2 = \begin{pmatrix} i & 0 \\ 0 & -i \end{pmatrix}$$

and
$$B_2 = \begin{pmatrix} i & -2i \\ 0 & -i \end{pmatrix}.$$

The fixed points of the generators are fix $A_2 = \{0, \infty\}$ and fix $A_2 = \{1, \infty\}$. The gluing construction is realized by adding a new generator $C_2[\tau]$ that conjugates the generators $A_2$ and $B_2$ and maps the fixed points as follows: $C_2[\tau](0) = \infty$, $C_2[\tau](\infty) = 1$:
$$C_2[\tau] = \begin{pmatrix} -\tau & 1/\tau \\ -\tau & 0 \end{pmatrix}.$$

It can be shown (as in [20] for hyperbolic triangle groups) that for large $|\tau| > 2$,
$$G_2[\tau] = \langle F_2, C_2[\tau] \rangle$$

is a discrete group with a connected set of discontinuity such that $\Omega(G_2[\tau])/G_2[\tau]$ is a punctured torus, and that the map $G_2[\tau] \mapsto \tau^2$ is a global complex analytic coordinate in the deformation space of the group $G_2[\tau_0]$ for any $\tau_0$ that gives a group uniformizing a punctured torus.

As in the case $n > 2$ we denote the image of the deformation space of $G_2[\tau_0]$ by $\mathcal{M}_2$.

*3.4. Remarks.* (1) The proof of Theorem 1 of [20] can be used to show that $\mathcal{M}_2$ is naturally isomorphic to the Riley slice $\mathcal{R}$ of Schottky space treated in [8].

(2) Figure 3 shows a $2-1$ regular covering of the deformation space $\mathcal{M}_{21}$: The plotted parameter is $i \operatorname{tr} C_\tau = i \sinh(d_n/2)(\tau - 1/\tau)$. $\mathcal{M}_n \subset \{\tau^2 : |\tau| > 1\}$, and the map $\tau \mapsto \tau - \frac{1}{\tau}$ maps the outside of the unit disk injectively onto the complement of the interval $[-2i, 2i]$.

$$\{\tau \in \mathbb{C} : |\tau| > 1\} \xrightarrow[1-1]{\tau \mapsto \tau - \frac{1}{\tau}} \mathbb{C} \setminus [-2\sinh(d_n/2), 2\sinh(d_n/2)] \supset \{\operatorname{tr} T[\tau] : \tau^2 \in \mathcal{M}_n\}$$

$$z \mapsto z^2 \downarrow 1-2$$

$$\{\tau : |\tau| > 1\} \supset \mathcal{M}_n$$

Thus, it is clear that the map $\operatorname{tr} C[\tau] \mapsto \tau^2$ is 2 to 1 from $\{\operatorname{tr} T[\tau] : \tau^2 \in \mathcal{M}_n\}$.

If $\tau^2 \in \mathcal{M}_n$, $n$ finite, the invariant component $\Omega_0$ of $G_n[\tau]$ is infinitely connected. There are homotopically nontrivial simple closed curves in $\Omega_0$ that are stabilized by conjugates of the finite cyclic group generated by $A_n$. These loops project $n$-to-1 to a non-dividing simple closed geodesic $\alpha$ on the punctured torus $\Omega_0/G_n[\tau]$. Thus the deformation space $\mathcal{T}(G_n[\tau])$ is $\mathcal{T}(1,1)$ factored by the action of the subgroup of the modular group generated by the $n$th power of the Dehn twist along this special curve. Complex analytically $\mathcal{T}(G_n[\tau])$, $2 \leq n < \infty$, is a punctured disk, in $\mathcal{M}_n$ the puncture is $\infty$. More detailed treatments of these deformation spaces can be found in [1] and [20]

## 4. Circle chains and rational pleating rays for deformation spaces of Koebe groups

In this section we sketch, following [7] and [8], the definitions and basic properties of rational pleating rays for $\mathcal{M}_n$, $3 \leq n < \infty$. The introduction of pleating structure allows us to state and prove Theorem 6.6 that identifies the groups in $\mathbb{H} \setminus \overline{\mathscr{M}}$ that are discrete by Theorem 6.5, as groups that appear in the boundary of these deformation spaces as the endpoints of properly chosen rational pleating rays. The corresponding structure for $\mathscr{M} = \mathscr{M}_\infty$ was given defined in [7]. Keen and Series treated the pleating structure of the Riley slice $\mathscr{R}$ in [8]. By Remark 3.4(1) we can use this to define the pleating structure on $\mathscr{M}_2$. For the purpose of proving Theorem 6.6 the case $n \geq 3$ is the most important, and we concentrate on that. The development follows [7] closely and we omit most of the details.

Let $G_n[\tau] = \langle A_n, C[\tau] \rangle$ be a Koebe group uniformizing a punctured torus. Let $W_{p/q}[A_n, C[\tau]]$ be the group element in $G_n[\tau]$ corresponding to the word in $A_n$ and $C[\tau]$ formed by the Farey procedure as in Section 1. Clearly some of the different words project to the same group element, that is, we do not get a 1-1 correspondence between the elements represented by the words $W_{p/q}$ and free homotopy classes of simple closed curves on the punctured torus. However, the identifications take place in a controlled manner, and we have

**4.1. Lemma.** $W_{p/q}[A_n, C[\tau]] = W_{p'/q'}[A_n, C[\tau]]$ *if and only if* $p'/q' = p/q + kn$ *for some* $k \in \mathbb{Z}$.

*Proof.* The $p/q$ words have the following nice property (See Osborne and Zieschang [19], Theorem 3.5):

$$W_{p/q+m}[X,Y] = W_{p/q}[W_m[X,Y], W_{m+1}[X,Y]].$$

Using this and the fact $W_m[X,Y] = X^{-1}Y$ shows that it is enough to prove that if $0 < p/q < p'/q' < 1$, then $W_{p/q}[A_n, C] \neq W_{p'/q'}[A_n, C]$, but this is clear, as $A_n$ only appears in the powers $0$ and $-1$ for these words, and the fact that $A_n^n = \mathrm{id}$ causes no reductions. □

For the rest of this section rational numbers $p/q$ should be thought of as points in $\mathbb{Q}/n\mathbb{Z}$. Lemma 4.1 implies that we get an enumeration by $\mathbb{Q}/n\mathbb{Z}$ of the set of equivalence classes of free homotopy classes of simple closed curves on the punctured torus, where the equivalence relation is generated by the relation $\alpha^n$, $\alpha$ being the simple closed curve on the torus corresponding to the elliptic element $A_n$. Accordingly, we call any simple closed curve in the class of the projection of $W_{p/q}$ a $p/q \mod n\mathbb{Z}$-*curve*, or a $p/q$-*curve*.

For each $p/q \in \mathbb{Q}/n\mathbb{Z}$ we define the $p/q$ *pleating ray* for the group $G_n[\tau]$ to be

$$\mathscr{P}_{p/q}^n = \left\{ \tau^2 \in \mathscr{M}_n : \partial \mathscr{C}_0(G^n[\tau^2])/G^n[\tau^2] \text{ is pleated along a } p/q \mod n\mathbb{Z}\text{-curve} \right\}.$$

As in the case of terminal regular b-groups it follows from [7] Lemma 4.6 that each rational ray is contained in the real locus of the trace of the special word $W_{p/q}$

$$\mathscr{P}_{p/q}^n \subset \widetilde{\mathscr{V}}_{p/q}^n = \left\{ \tau^2 \in \mathscr{M}_n : \operatorname{tr} W_{p/q} \in \mathbb{R} \right\}.$$

If a class contains a curve corresponding to $W_{p/q}$, for which $p/q \in \mathbb{Z}$, all curves in the class correspond to integral values of $p/q$. Thus, we call $\mathscr{P}_{m/1}^n$ an *integral pleating ray*. We can generalize the concept of circle chains used by Keen–Series and Wright:

**4.2. Definition.** Let $G = \langle F_n, C \rangle$ be a group of Möbius transformations (not necessarily discrete), where $F_n = \langle A_n, B_n \rangle$ as earlier, and $C$ is a transformation satisfying
$$C^{-1}A_n C = B_n^{-1}.$$

Let $W_{p/q} = W_{p/q}[A, C]$. A collection $\{\delta_i\}$, $i \in \mathbb{Z}$ mod $nq$ of closed, round disks is a *combinatorial $p/q$ chain* if it satisfies the following conditions:

(1) $\delta_0$ is tangent to $\Lambda(F)$, the limit set of $F$, at the fixed point of the parabolic element
$$K = W_{p/q} W_{r/s}^{-1} W_{p/q}^{-1} W_{r/s},$$

(2) $W_{p/q}\delta_0 = \delta_0$,
(3) $C(\delta_j) = \delta_{j+p}$ for all $j = 0, \ldots, q$,
(4) $A(\delta_j) = \delta_{j+q}$ for all $j = 0, \ldots, p$, and
(5) $\delta_i \neq \delta_j$ if $j \neq i$.

We say the chain *connects* fix $K$ to $A(\text{fix } K)$. The chain is *tangent* if the If the interiors of the disks are pairwise disjoint but $\delta_i$ and $\delta_{i+1}$ are tangent. If we know that the group is discrete, we say that the chain is *proper* if the interiors of the disks $\delta_i$ are contained in $\Omega(G)$ for all $i$, adjacent disks $\delta_i$ and $\delta_{i+1}$ intersect for all $i$, and $\text{int } \delta_i \cap \text{int } \delta_j = \emptyset$ for $|i - j| > 1$.

*4.3. Remarks.* (1) The combinatorial $p/q$ chain is associated to a group with a fixed pair of generators. For a different choice of generators there may be a chain with different combinatorics, this is used in the proofs of Theorems 4.11 and 6.6, see also McShane, Parker and Redfern [18].

(2) Note that we do not require that the triangle group be hyperbolic. For $n = 2$ the special word $W_{p/q}$ is an elliptic of order 2, and $F_2$ is the infinite dihedral group.

(3) The disks $\delta_0, \ldots \delta_{p+q-1}$ form a combinatorially interesting set: This part can be generated starting from $\delta_0$ by the following process:

(i) If $0 \leq i < q$, set $\delta_{i+p} = T(\delta_i)$,
(ii) If $q \leq i < p + q$, set $\delta_{i-q} = S^{-1}(\delta_i)$.

(4) Unlike Keen and Series [7] we do not require that the group $G$ is discrete. In fact, we use the existence of special circle chains to prove that the group $G[\mu]$ defined in Section 1 is discrete for a number of the parameters $\mu$ found in Lemma 1.3.

In [7] proper chains were used to analyze the structure of the invariant sets of terminal b-groups and the geometry of their convex hulls. Wright [25] used tangent chains to construct fundamental sets for groups on the boundary of $\mathcal{M}$.

Using the methods of [7] we get the following properties of circle chains and rational pleating rays in $\mathcal{M}_n$:

**4.4. Lemma.** *The group $G_n[\tau]$ has a proper $p/q$ chain if and only if $\tau^2 \in \mathscr{P}_{p/q}^n$.*

*Proof.* [7] Proposition 4.11.

**4.5. Lemma.** *Let $0 \leq m < n$. The integral pleating ray $\mathscr{P}_{m/1}^n$ is the radial line*
$$\mathscr{P}_{m/1}^n = \mathcal{M}_n \cap \left\{ \tau^2 = t^2 e^{-i2\pi m/n} : t > 0 \right\}$$

in $\mathscr{T}(G_n[\tau])$.

*Proof.* It is easy to check that

$$W_{m/1} = \begin{pmatrix} e^{i\pi m/n}\sinh(d/2)\tau & * \\ * & -e^{-i\pi m/n}\sinh(d/2)/\tau \end{pmatrix},$$

so if we write $\tau = te^{i\varphi}$, we get

$$\operatorname{Im}\operatorname{tr} W_{m/1} = \left(t + \frac{1}{t}\right)\sinh(d/2)\sin(\varphi + \pi m/n).$$

Thus

$$\widetilde{\mathscr{V}}^n_{m/1} \subset \{\tau^2 \in \mathscr{M}_n : \arg(\tau^2) = -2\pi m/n\}.$$

It remains to show the existence of an $m/1$ chain. Here we can follow the argument of [7] Lemma 5.2. □

For the real locus of the trace of the general $p/q$-word we have (as in [7] Lemma 5.3)

**4.6. Lemma.** *Let $0 < p/q < n$, $n \notin \mathbb{Z}$. In the sector*

$$\{\tau^2 \in \mathscr{M}_n : -2\pi(k+1)/n \leq \arg\tau^2 \leq -2\pi k/n\}$$

*the real locus $\widetilde{\mathscr{V}}^n_{p/q}$ is asymptotic to the line*

$$\{\tau^2 \in \mathscr{M}_n : \arg\tau^2 = -2\pi p/(qn)\}.$$

*Proof.* An easy induction argument as in [7] Proposition 3.1 shows that

$$W_{p/q} = \begin{pmatrix} e^{p\,i\pi/n}\tau^q\sinh(d/2)^q + O(|\tau^{q-1}|) & O(|\tau|^q) \\ O(|\tau|^{q-2}) & O(|\tau|^{q-2}) \end{pmatrix}$$

as $|\tau| \to \infty$. When $|\tau| \to \infty$ the condition $\operatorname{tr} W_{p/q} \in \mathbb{R}$ approaches

$$\arg\tau^2 = -2\pi p/(nq). \quad \square$$

To show that the pleating ray is in the correct sector in $\mathscr{M}_n$ we need the continuity of the pleating map [9] Theorem 4. We denote by $\mathscr{PML}$ the space of projected measured laminations on the once-punctured torus. This space is the completion of the space of simple closed geodesics in the punctured torus with the topology induced from the embedding in $\widehat{\mathbb{R}}$ defined by the words $W_{p/q}$, thus $\mathscr{PML}$ is canonically identified with $\widehat{\mathbb{R}} = \mathbb{S}^1$.

**4.7. Lemma.** *Let $0 < p/q < n$, $n \notin \mathbb{Z}$. Then $\mathscr{P}_{p/q}^n$, $n \geq 3$, is contained in*

$$\mathscr{V}_{p/q}^n = \left\{\tau^2 \in \widetilde{\mathscr{V}}_{p/q}^n \cap \mathscr{M}_n : -2\pi(k+1)/n \leq \arg \tau^2 \leq -2\pi k/n\right\},$$

*where $k$ is the integral part of $p/q \mod n$.*

*Proof.* The map $\mathrm{pl} \colon \mathscr{M} \to \mathscr{PML} = \widehat{\mathbb{R}}$ that maps the parameter $\tau^2$ in the deformation space to the pleating locus of $\partial \mathscr{C}_0$ is continuous, see Keen and Series [9] Theorem 4, and [7] Theorem 4.5. Also, $\mathrm{pl}(\mathscr{M}) = \mathbb{R} \subset \widehat{\mathbb{R}}$.

The Dehn twist along the geodesic $\alpha$ corresponding to the element $S \in G[\mu]$, denoted by $\mathrm{Dehn}\,\alpha$, acts on both $\mathscr{M}$ and $\mathscr{PML}$, and the actions are compatible:

$$\mathrm{Dehn}_\alpha(\mu) = \mu + 2, \quad \mathrm{Dehn}_\alpha(r) = r + 1$$

for $\mu \in \mathscr{M}$, $r \in \mathscr{PML} = \widehat{\mathbb{R}}$. Note that $\mathrm{Dehn}_\alpha$ fixes $\infty \in \mathscr{PML}$ and acts freely discontinuously in $\mathrm{pl}(\mathscr{M}) = \mathbb{R}$. Thus the pleating map $\mathrm{pl}_n \colon \mathscr{M} \to \mathbb{S}^1$ becomes well defined and continuous:

$$\begin{array}{ccc} \mathscr{M} & \xrightarrow{\mathrm{pl}} & \mathbb{R} \subset \widehat{\mathbb{R}} = \mathscr{PML} \\ \text{projection} \downarrow & & \downarrow \text{projection} \\ \mathscr{M}_n & \xrightarrow{\mathrm{pl}_n} & \mathbb{S}^1 \end{array}$$

By lemma 4.5 we know that the integral pleating rays $\mathscr{P}_m^n \subset \mathscr{M}_n$ are exactly the radial lines of argument $-i2\pi m/n$.

The claim follows from these facts by considering the restriction of $\mathrm{pl}_n$ to a circle with radius big enough so that it is contained in $\mathscr{M}_n$: Continuity implies that $\mathscr{P}_{p/q}^n$ is non-empty and it has to be in the unique sector between $\mathscr{P}_k^n$ and $\mathscr{P}_{k+1}^n$ containing no integral rays. $\square$

Combining these observations with the analysis of Keen and Series [8], [7] Section 5 we get

**4.8. Theorem.** *For $p/q \in \mathbb{Q}/n\mathbb{Z}$ the pleating ray $\mathscr{P}_{p/q}^n$ coincides with the branch of $\widetilde{\mathscr{V}}_{p/q}^n$ with asymptotic argument $\arg \tau^2 = -2\pi i p/(qn)$. This branch contains no singularities and $\overline{\mathscr{P}}_{p/q}^n \cap \partial \mathscr{M}_n$ consists of a single point $\tau_{p/q}$ such that the transformation $W_{p/q}[A_n, C_n[\tau_{p/q}]]$ is parabolic. Furthermore, the group $G_n[\tau_{p/q}]$ is geometrically finite, and it has a tangent circle chain connecting $\mathrm{fix}\, K_n$ to $A_n \mathrm{fix}\, K_n$.*

*4.9. Remark.* One can define extended rational pleating rays for $\mathscr{M}_n$ as was done in Section 1 for $\mathscr{M}$ and ask if there is a result analogous to Theorem 6.5 for these rays. However, we do not pursue this, because the method we use in Sections 5 and 6 to treat discreteness of the groups on the extensions cannot be used; The method relies on the fact that $\langle S, \widetilde{S} \rangle$ is a subgroup of $G[\mu]$ for any parameter $\mu \in \mathbb{C}$.

We now prove a special case of our main Theorems: The integral rays, that is the pleating rays $\mathscr{P}_{m/1} = \mathscr{P}_m$, $m \in \mathbb{Z}$ are relatively easy to handle, they are the extended $m/1$ ray is the line $\{z : \mathrm{Im}\, z = 2m\}$, and the only possible values of $\mu$ on these rays for which $G[\mu]$ is discrete are clearly of the form

$$\mu_m(n) = 2m + 2i\cos(\pi/n).$$

The following theorem shows that, in fact, all these parameters give discrete groups and describes the quotient Riemann surfaces of these groups. This result will be generalized in the final section, where we treat all rational pleating rays. However, for non-integral rays the result is much more difficult to prove by our method, and we only get a local result. The following sections develop the machinery required for non-integral rays.

**4.10. Theorem.** *Let $\mu_m(n) = 2m + 2i\cos(\pi/n)$. Then*
  (1) *$G[\mu_m(n)]$ is discrete and geometrically finite,*
  (2) *$G[\mu_m(n)] = \langle T[\mu_m(n)], S^{-1}T[\mu_m(n)]S\rangle *_S$,*
  (3) *if $n \geq 2$, then $\Omega(G[\mu_m(n)])/G[\mu_m(n)]$ consists of a punctured sphere and a sphere with one puncture and two elliptic points of order $n$,*
  (4) *$\Omega(G[\mu_m(2)])/G[\mu_m(2)]$ is a thrice punctured sphere.*

*Proof.* The groups $G[\mu_m(n)]$, $m \in \mathbb{Z}$ are all isomorphic to $G[\mu_0(n)]$. Thus, it is enough to consider the parameters $\mu_0(n) = 2i\cos(\pi/n)$. Let $T = T[2i\cos(\pi/n)]$. $T$ is an elliptic of order $n$. By Lemma 2.2 we see that the group

$$F = \langle T, S^{-1}T^{-1}S\rangle$$

is a triangle group of signature $(n, n, \infty)$.

*Case $n > 2$:* The set

$$\mathscr{F} = \widehat{\mathbb{C}}\backslash\big(\{z : |z| \leq 1\} \cup \{z : |z + 2i\cos(\pi/n)| \leq 1\}$$
$$\cup \{z : |z + 2| \leq 1\} \cup \{z : |z + 2 + 2i\cos(\pi/n)| \leq 1\}\big)$$

is a fundamental polygon for the action of $F$ on

$$\Omega(F) = \widehat{\mathbb{C}} \setminus \{z : |z + 2 + i\cos(\pi/n)\}.$$

Let

$$\mathscr{D}_1 = \{z : \operatorname{Re} z < -2\}, \quad \mathscr{D}_2 = \{z : \operatorname{Re} z > 0\}.$$

It is easy to see that $T$ preserves $\mathscr{D}_2$. Also, since $\mathscr{F} \cap \mathscr{D}_2$ is a fundamental polygon for the action of $\langle T\rangle$ in $\mathscr{D}_2$, $\mathscr{D}_2$ is precisely invariant under $\langle T\rangle$ by Maskit [16] II.I.3. Also, $S(\operatorname{int}\mathscr{D}_1) = \operatorname{ext}\mathscr{D}_2$. We will now apply Maskit's second combination theorem (Version IV), [17]. The assumptions of the combination theorem are satisfied: $T(\operatorname{int}D_1) = \operatorname{ext}\mathscr{D}_2$, and for the sets $\Theta_1$ and $\Theta_2$ in the combination theorem we can clearly choose

$$\Theta_1 = \{S^{-1}T^iS(\infty) : i = 1, 2, \ldots n\},$$

and

$$\Theta_2 = \{T^i(\infty) : i = 1, 2, \ldots n\} = S(\Theta_1).$$

The statements now follow directly from the second combination theorem.

*Case $n = 2$:* The argument of the case $n > 2$ applies with the following modifications: Now $T(z) = 1/z$, and the group $F$ is the infinite dihedral group. The set

$$\mathscr{F} = \widehat{\mathbb{C}} \setminus \big(\{z : |z| \leq 1\} \cup \{z : |z + 2| \leq 1\}\big)$$

is a fundamental polygon for $F$ acting on $\widehat{\mathbb{C}} \setminus \{-1\}$. $\square$

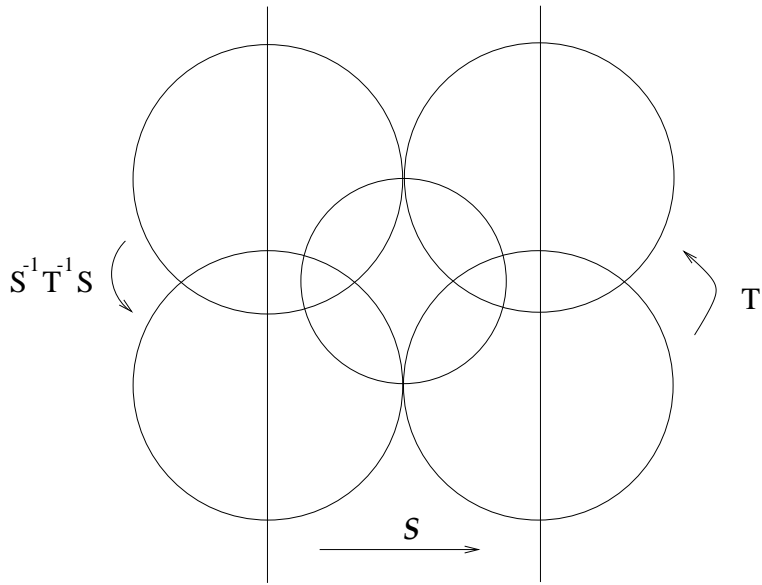

FIGURE 4. The construction of the fundamental set in Theorem 4.10 for $n \geq 3$.

**4.11. Theorem.** *Let $\tau_m$ be the unique boundary point of $\mathcal{M}$ on $\mathscr{P}_m^n$. Then there is a Möbius transformation $\beta$ such that $\beta G_n[\tau_m]\beta^{-1} = G[\mu_m(n)]$*

*Proof.* Let $\beta$ be the Möbius transformation that satisfies
$$\beta(ie^{-i\pi/n} - 2) = 0, \quad \beta(ie^{i\pi/n} - 2) = \infty, \quad \text{and} \quad \beta(-1) = e^{-i\pi/n}.$$

Now it is a routine calculation to check that
$$\beta(S^{-1}T^{-1}S)\beta^{-1} = A_n,$$
$$\beta T \beta^{-1} = B_n,$$
$$\beta S \beta^{-1} = C_n[\tau_{0/1}],$$

where
$$\tau_{0/1}^2 = \frac{1 + \sin(\pi/n)}{1 - \sin(\pi/n)}$$

is the boundary point of the 0/1 pleating ray $\mathscr{P}_{0/1}^n$ in $\mathscr{M}_n$. □

## 5. TANGENT CIRCLE CHAINS AND DISCRETENESS

In this Section we prove theorem 5.1, a technical result used in the proof of Theorem 6.5. We show how to use the existence of a tangent circle chain and Maskit's second combination theorem [17] to establish the discreteness of a group of the form $G = \langle F_n, C \rangle$. Clearly, there can be non-discrete groups that possess tangent circle chains. We will work with tangent circle chains that satisfy the additional property

(∗) $\qquad\qquad\qquad \delta_i \cap \delta_j = \emptyset \quad \text{if} \quad |i - j| > 1.$

In Proposition 6.3 we show that the circle chains of cusp groups $G[\mu_{p/q}]$ satisfy (∗). In the proof of Theorem 6.5 we produce circle chains satisfying (∗) for parameters on $\mathscr{P}_{p/q}^+$ by perturbing the chains of $G[\mu_{p/q}]$.

Let $\delta_i$, $i \in \mathbb{Z} \mod n$ be a tangent $p/q$ circle chain for $G$ satisfying $(*)$. Let $\gamma_i$ be the unique circular arc orthogonal to $\partial \delta_i$ connecting the points of tangency of $\delta_i$ with $\delta_{i-1}$ and $\delta_{i-1}$, $i \in \mathbb{Z} \mod n$. Let

$$\mathscr{W}_A = \bigcup_{i \in \mathbb{Z} \mod n} \gamma_i.$$

Let

(†) $$\mathscr{W}_B = C^{-1}(\mathscr{W}_A).$$

Condition $(*)$ implies that $\mathscr{W}_A$ and $\mathscr{W}_B$ are Jordan curves. Let $\mathscr{D}_A$ be the component of $\widehat{\mathbb{C}} \setminus \mathscr{W}_A$ contained in the same component of $\widehat{\mathbb{C}} \setminus \Lambda(F)$ as $\mathscr{W}_A$. Let $\mathscr{D}_B$ be the component of $\widehat{\mathbb{C}} \setminus \mathscr{W}_B$ contained in the same component of $\widehat{\mathbb{C}} \setminus \Lambda(F)$ as $\mathscr{W}_B$.

**5.1. Theorem.** *Let $G = \langle F, C \rangle$, where $F = \langle A, B \rangle$ is a triangle group of signature $(n, n, \infty)$, and $C$ a transformation satisfying*

$$C^{-1}AC = B^{-1}.$$

*If $G$ has a tangent $p/q$- chain satisfying $(*)$ such that*

(i) *the interiors of the disks $\mathscr{D}_A$ and $\mathscr{D}_B$ are disjoint, and*
(ii) *the transformations $W_{p/q}[A, C]$ are parabolic,*

*then*

(1) $G = F *_C$,
(2) $G$ *is discrete,*
(3) $G$ *is geometrically finite,*
(4) *If $n > 3$, $\Omega(G)/G = S_1 \cup S_2$, where $S_1$ has signature $(0, 3; n, n, \infty)$, and $S_2$ has signature $(0, 3; \infty, \infty, \infty)$. If $n = 2$, then $\Omega(G)/G = S_2$.*

The proof of this theorem is an application of Maskit's second combination theorem [17]. The following definition is from Maskit [16]:

**5.2. Definition.** Let $\Gamma$ be a Kleinian group, $J_1$, $J_2$ subgroups, and $C$ a loxodromic transformation such that $CJ_1C^{-1}$. Three disjoint nonempty sets $(X_1, X_2, Z)$ form an *interactive triple* for $(\Gamma, C)$ (or $(\Gamma, J_1, J_2, C)$), if

(1) $X_1$ is precisely invariant in $\Gamma$ under $\langle J_1 \rangle$,
(2) $X_2$ is precisely invariant in $\Gamma$ under $\langle J_2 \rangle$,
(3) $C(Z \cup X_2) \subset X_2$ and $C^{-1}(Z \cup X_1) \subset X_1$,
(4) For all $g \in \Gamma$ $g(X_1) \subset Z \cup X_1$, and $g(X_2) \subset Z \cup X_2$.

The interactive triple $(X_1, X_2, Z)$ is *proper* if $Z \setminus \Gamma(X_1 \cup X_2) \neq \emptyset$.

**5.3. Lemma.** *Let $\mathscr{D}_A$, $\mathscr{D}_B$, $F$ and $C$ be as in Theorem 5.1. Then*

$$(\operatorname{int} \mathscr{D}_A, \operatorname{int} \mathscr{D}_B, \widehat{\mathbb{C}} \setminus (\mathscr{D}_A \cup \mathscr{D}_B))$$

*is an interactive triple for $(F, C)$.*

*Proof.* Property (3) is true by assumption of the Theorem. Also, it is clear that $\mathscr{D}_A$ is $A$-invariant, and that $\mathscr{D}_B$ is $B$-invariant. By Maskit [16] Propositions VII.E.3 and VII.E.4 it is enough to show that $\mathscr{D}_A$ is precisely $\langle A \rangle$ invariant and that $\mathscr{D}_B$ is precisely $\langle B \rangle$ invariant in $F$, that is, for all $g \in F \setminus \langle A \rangle$ $g(\operatorname{int} \mathscr{D}_A) \cap \operatorname{int} \mathscr{D}_A = \emptyset$, and similarly for $\mathscr{D}_B$. By Maskit [16] Proposition II.I.3 it is enough to show that there is a fundamental set $P$ for $F$ such that $P \cap \mathscr{D}_A$ is a fundamental set for the action of $\langle A \rangle$ in int $\mathscr{D}_A$.

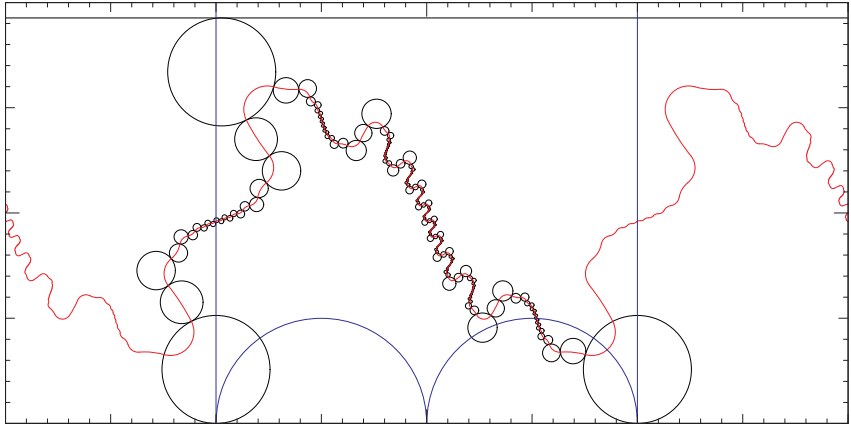

FIGURE 5. The circle chain for $G[\mu_{21/229}]$. The curve $\mathscr{W}_A$ intersects the disk $\{|z - \frac{1}{2}| < \frac{1}{2}\}$. One obtains the fundamental set $P$ by adding the set $\mathscr{D}_A \cap \{|z - \frac{1}{2}| < \frac{1}{2}\}$ to the standard fundamental polygon $P_0 = \{z \in \mathbb{C} : -1 \leq \operatorname{Im} z < 1, |z + \frac{1}{2}| \geq \frac{1}{2}, |z - \frac{1}{2}| > \frac{1}{2}\}$ and by removing its image $B\left(\mathscr{D}_A \cap \{|z - \frac{1}{2}| < \frac{1}{2}\}\right)$.

We use the fact that $B = C^{-1}A^{-1}C$: By assumption $C(\mathscr{D}_A) \subset \mathscr{D}_A$. By $A$ invariance $A^{-1}C(\mathscr{D}_A) \subset \mathscr{D}_A$. Now $A^{-1}C(\delta_{q-1}) = \delta_{p-1}$. Thus, $C^{-1}A^{-1}C(\delta_{q-1})$ meets $\mathscr{W}_B$, but not $\mathscr{W}_A$. Similarly, we see that for all $i = 1, \ldots, q-1$ the part of the boundary of $\mathscr{D}_A$ inside $\delta_i$ is mapped by $B$ into the part of the boundary of $\mathscr{W}_B$ not touching $\mathscr{W}_A$. For $B^{-1}$ we have the corresponding result by the same kind of reasoning: For $i = 1, \ldots, q-1$, we have $AC(\delta_i) = \delta_{p+q+i}$. As $i > 0$, the combinatorics of the $p/q$ chain gives that the part of $\mathscr{W}_A$ inside $\delta_i$ is mapped into the part of $\mathscr{W}_B$ not touching $\mathscr{W}_A$. Thus we have that the images $B^{\pm 1}(\operatorname{int} \mathscr{D}_A)$ are disjoint from $\mathscr{D}_A$. We can now modify the 'standard' fundamental polygon of the triangle group $F$ to obtain the required fundamental set (See Figure 5). □

Next we will modify the Jordan curves so that they meet the requirements of the fourth version of Maskit's second combination theorem [17]: $\mathscr{W}_A$ and $\mathscr{W}_B$ are only allowed to meet at the points of tangency of the disks $\delta_0, \ldots, \delta_q$, not along a union of circular arcs as above. If $p/q \in \mathbb{Z}$, $\mathscr{W}_A$ and $\mathscr{W}_B$ are round disks and we do not change them. Let $p/q \in \mathbb{Q} \setminus \mathbb{Z}$. As groups without reference to a specific pair of generators $G[\mu] = G[\mu + 2i]$, and we can thus assume $0 < p/q < 1$.

Let $\gamma$ be a directed geodesic in $\mathbb{D}$ connecting the points $z_1$ and $z_2$ on $\partial \mathbb{D}$, $z$ a point on $\gamma$, and $\epsilon > 0$. Let $z(\epsilon)$ be the point on the geodesic perpendicular to $\gamma$ at $z$ that lies the distance $\epsilon$ to the right from $z$. We denote by $\gamma(\epsilon)$ the union of the geodesic arcs connecting $z_1$ to $z(\epsilon)$ and $z(\epsilon)$ to $z_2$.

Enumerate the arcs $\widetilde{\gamma}_i$ that make up $\mathscr{W}_B$ so that $\widetilde{\gamma}_i \subset \delta_i$ for $i = 1, \ldots q-1$. In the modification process we go through the disks $\delta_0, \ldots, \delta_q$ respecting the combinatorics described in Remark 4.3(3), we replace $\mathscr{W}_A$ by $\mathscr{W}_A'$ and $\mathscr{W}_B$ by $\mathscr{W}_B'$ as follows: Let

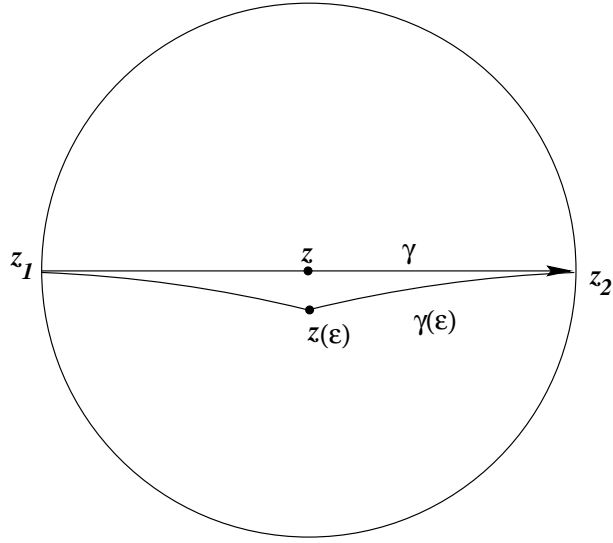

FIGURE 6. The modification of $\mathscr{W}_A$ inside a disk of the circle chain.

$\epsilon > 0$. Think of each $\delta_i$ as the hyperbolic plane, and orient the geodesic $\gamma_i$ so that it starts at $\delta_{i-1} \cap \delta_i$ and ends at $\delta_i \cap \delta_{i+1}$. First set $\widetilde{\gamma}'_0 = \gamma_0$, $\gamma'_p = C(\widetilde{\gamma}'_p) = \gamma_p$. As the first modification set $\widetilde{\gamma}'_p = \gamma_p(\epsilon)$. Go through the chain as follows: Assume you have made the modification inside $\delta_i$, and that $\widetilde{\gamma}'_i = \gamma_p(k\epsilon)$. There are three cases:

(1) If $i + p < q$ set $\gamma'_{p+i} = C(\widetilde{\gamma}'_i)$, $\widetilde{\gamma}'_{i+p} = \gamma_{i+p}((k+1)\epsilon)$. Continue modification.
(2) If $i + p > q$ set $\gamma'_{p+i-q} = A^{-1}C(\widetilde{\gamma}'_i)$, $\widetilde{\gamma}'_{i+p-q} = \gamma_{i+p}((k+1)\epsilon)$. Continue modification.
(3) If $i + p = q$ set $\widetilde{\gamma}'_q = \gamma_p((k+1)\epsilon)$. Now we have changed

After these modifications we set

$$\mathscr{W}'_A = \bigcup_{n \in \mathbb{Z}} \bigcup_{i=1}^{q} A^n(\gamma'_i), \quad \mathscr{W}'_B = \bigcup_{n \in \mathbb{Z}} \bigcup_{i=0}^{q-1} B^n(\widetilde{\gamma}'_i).$$

It is clear from the construction that we can choose $\epsilon$ in such a way that the curves $\mathscr{W}'_A$ and $\mathscr{W}'_B$ satisfy the same mapping properties with respect to $G$ and $C$ as the original curves $\mathscr{W}_A$ and $\mathscr{W}_B$. To complete the modification we replace $\mathscr{D}_A$ and $\mathscr{D}_B$ by $\mathscr{D}'_A$ and $\mathscr{D}'_B$ such that $\partial \mathscr{D}'_A = \mathscr{W}'_A$ and $\partial \mathscr{D}'_B = \mathscr{W}'_B$. Clearly, as in Lemma, $\left(\operatorname{int} \mathscr{D}'_A, \operatorname{int} \mathscr{D}'_B, \widehat{\mathbb{C}} \setminus (\mathscr{D}'_A \cup \mathscr{D}'_A)\right)$ is an interactive triple.

*Proof of Theorem 5.1.* Claims $(1) - (3)$: In the notation of Maskit's second combination theorem (Version IV) [17], Section II,

$$G_0 = F, \quad \langle B \rangle = J_1, \quad \langle A \rangle = J_2,$$
$$f = C, \quad B_1 = \mathscr{D}'_B, \quad \text{and } B_2 = \mathscr{D}'_A.$$

Set
$$\Theta_A = \bigcup_{i \in \mathbb{Z}} A^i(\mathscr{W}'_A \cup \mathscr{W}'_B), \quad \Theta_B = \bigcup_{i \in \mathbb{Z}} B^i(\mathscr{W}'_A \cup \mathscr{W}'_B).$$

By the assumptions of the Theorem, Lemma 5.3, and the modification above, the assumptions of the second combination theorem are satisfied:

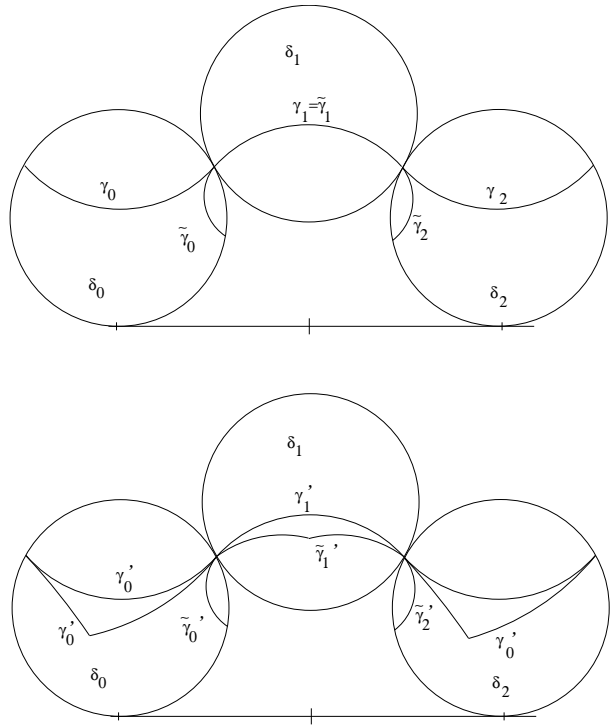

FIGURE 7. An example of the modification of $\mathscr{W}_A$ and $\mathscr{W}_B$ in the case $p/q = 1/2$.

Assumption (A) is satisfied, except that the curves $\mathscr{W}'_A$ and $\mathscr{W}'_B$ are not locally circular at $\Theta_A$ and $\Theta_B$, respectively. However, as remarked in [17] Section 0.2.3, this requirement can be weakened considerably, in our situation, see Figures 6 and 7, the curves consist of circular arcs that are perpendicular to the boundaries of the disks $\delta_i$ at the points of $\Theta_A$ and $\Theta_B$. Also, we can make a small local perturbation of the curves $\mathscr{W}'_A$ and $\mathscr{W}'_B$ at their points of tangency such that the curves are locally circular at these points without introducing new points of tangency, or changing the mapping properties ((†) and $(1) - (4)$ of Definition 5.2)

Assumption (B) of the combination theorem is exactly assumption (i) of Theorem 5.1: The cyclic stabilizers of Maskit [17] Section II.2 are conjugate to $W_{p/q}$. Assumption (C), $\widehat{\mathbb{C}} \setminus (\mathscr{D}_A \cup \mathscr{D}_B)$ is trivially satisfied.

Now, claim (1) is Maskit's conclusion (i), claim (2) is (ii), and claim (3) is (x).

*Proof of claim* (4): If $G$ is a Kleinian group such that the limit set $\Lambda(G)$ consists of more than one point (like F in our situation), we follow Maskit [16] and denote by area$(G)$ the area in the hyperbolic metric of $\Omega(G)/G$. If $\Lambda(G)$ is finite, we set area$(G) = 0$, and if $G$ is a finite Kleinian group, area$(G) = -4\pi/\#G$. Conclusion (xiii) of Maskit's second combination theorem [17] gives in our situation for $n > 3$ the following:

$$\text{area}(G) = \text{area}(F) - \text{area}(\langle A \rangle) = 2 \cdot 2\pi \left(1 - \frac{2}{n}\right) + 4\pi/n = 2\pi \left(1 - \frac{2}{n}\right) + 2\pi.$$

By (1) we know that the parabolics $A$, $B$, and $W_{p/q}$ define different conjugacy classes. By Maskit's conclusion (viii) we know that the quotient contains one sphere with a puncture and two cone points of order $n$, and that $A$, $B$, $W_{p/q}$ correspond to three punctures on the remainder of the quotient. It follows by considering the Euler characteristic that the quotient is as in claim (4).

Similarly, for $n = 2$ we have:

$$\text{area}(G) = \text{area}(F) - \text{area}(\langle A \rangle) = 4\pi/2 = 2\pi.$$

and the quotient is a thrice punctured sphere. $\square$

## 6. Proof of main results

Let $\mu_{p/q}$ be the cusp on $\overline{\mathscr{P}_{p/q}}$. The quotient $\Omega(G[\mu_{p/q}])/G[\mu_{p/q}]$ is the union of two thrice punctured spheres, one corresponding to $\mathbb{H}^*/\langle S, \widetilde{S} \rangle$, the other to $\delta_0/\langle W_{p/q}, W_{r/s}^{-1} W_{p/q}^{-1} W_{r/s} \rangle$. The main idea comes from a property of the cusps of $\mathscr{M}$; McShane, Parker and Redfern [18], Proposition 7.1, made the following useful observation:

**6.1. Proposition.** *$G[\mu_{p/q}]$ is the HNN extension of the torsion-free triangle group $\langle W_{p/q}, W_{r/s}^{-1} W_{p/q}^{-1} W_{r/s} \rangle$ by $W_{r/s}$,*

$$G[\mu_{p/q}] = \langle W_{p/q}, W_{r/s}^{-1} W_{p/q}^{-1} W_{r/s} \rangle *_{W_{r/s}}.$$

*Moreover,*

(1) *$G[\mu_{p/q}]$ is the $-s/q$ cusp in the deformation space of terminal b-groups of the form $\langle W_{p/q}, W_{r/s}^{-1} W_{p/q}^{-1} W_{r/s} \rangle *_Z$ representing a punctured torus in their invariant component.*
(2) *there is a tangent $-s/q$ chain $(\Delta_i)_{i \in \mathbb{Z}}$, $\text{int } \Delta_i \subset \Omega(G[\mu_{p/q}])$ for these generators, with $\Delta_0 \cap \delta_0 = \text{fix } \widetilde{K}$, $\Delta_q \cap \delta_0 = \text{fix } K$.*

We need the following Lemma in the Proof of Theorem 6.5:

**6.2 Proposition.** *For any $\mu \in \mathbb{H}$*

$$W_{-s/q}[W_{p/q}[S,T], W_{r/s}[S,T]] = W_{1/0}[S,T] = S^{-1}.$$

*Proof.* Osborne and Zieschang [19], Theorem 3.5, proved the following: For any $m/n$, and any neighbors $p/q$ and $r/s$ the following formula holds:

$$W_{m/n}[W_{p/q}[S,T], W_{r/s}[S,T]] = W_{(mr+np)/(ms+nq)}[S,T].$$

The claim is a special case of this formula. $\square$

For the remainder of this section, let us use the notation

$$X = W_{p/q}, \quad Y = W_{r/s}^{-1} W_{p/q}^{-1} W_{r/s}.$$

$\delta_0$ is the invariant disk of $F_{p/q} = \langle X, Y \rangle$ such that $\text{int } \delta_0$ is a component of $\Omega(G[\mu_{p/q}])$. Note also that $\mathbb{H}^* = \Delta_q$, and the stabilizer $\Gamma_i$ of $\Delta_i$ is conjugate in $G[\mu_{p/q}]$ to $\langle S, \widetilde{S} \rangle$ for all $i$.

The map $\mu \mapsto G[\mu]$ induces a homomorphism from $G[\mu_{p/q}]$ onto any $G[\mu]$, $\mu \in \mathbb{C}$. Let us denote by $\Gamma_i[\mu]$ the image of $\Gamma_i[\mu_{p/q}]$ by this homomorphism. As $\langle S, \widetilde{S} \rangle$ does not depend on the parameter $\mu$, in particular this means that the groups $\Gamma_i[\mu]$ are discrete and there is a naturally defined disk $\Delta_i[\mu]$ for each $\mu \in \mathbb{H}$, bounded by the limit set of $\Gamma_i[\mu]$. In particular, $\Delta_i[\mu_{p/q}] = \Delta_i$.

**6.3. Proposition.** (1): $\Delta_i[\mu]$ and $\Delta_{i+1}[\mu]$ are tangent for any $\mu \in \mathbb{H}$.
(2): If $|i - j| > 1$, then $\Delta_i[\mu_{p/q}] \cap \Delta_j[\mu_{p/q}] = \emptyset$.

*Proof.* (1): $\Delta_i$ and $\Delta_{i+1}$ are stabilized by the same parabolic transformation conjugate to $W_{p/q}$. Thus these disks are clearly tangent at the fixed point of this element.

(2): Assume $j > i + 1$ and that there is a point $x \in \Delta_i[\mu_{p/q}] \cap \Delta_j[\mu_{p/q}]$. By Maskit [15] Theorem 3:
$$\Lambda(\Gamma_i) \cap \Lambda(\Gamma_j) = \Lambda(\Gamma_i \cap \Gamma_j).$$

$\Delta_{[\mu_{p/q}]i}$ and $\Delta_j[\mu_{p/q}]$ are round disk components of $G[\mu_{p/q}]$, so they can intersect in exactly one point. Thus $\Gamma_i \cap \Gamma_j = \langle P \rangle$, where $P$ is a parabolic element fixing $x$. Thus it has to be conjugate to one of $X$, $Y$, $XY$ or $W_{-s/q}[X, Z]$. The three first parabolics in this list correspond to punctures on $\delta_0/F_{p/q}$, so $\Delta_i[\mu_{p/q}]$ and $\Delta_j[\mu_{p/q}]$ cannot be tangent at the fixed point of an element conjugate to any of these. Assume they are tangent at the fixed point of an element $g$ conjugate to $W_{p/q}$. For parameters close to $\mu_{p/q}$ $g$ is loxodromic. However, $G[\mu]$ and $G[\mu_{p/q}]$ are naturally isomorphic, and $\Gamma_i$ is conjugate to $\Gamma_q = \langle S, \widetilde{S} \rangle \ni P$, which does not depend on $\mu$. Thus $P$ is parabolic, a contradiction. $\square$

We need to show that for $\mu$ close to the cusp the disks $\Delta_0[\mu], \ldots, \Delta_q[\mu]$ satisfy (*). This will be a consequence of the following Lemma:

**6.4. Lemma.** *For $\epsilon > 0$ there is a $\eta > 0$ such that if $|\mu - \mu_0| < \eta$, then*
$$\max\{z \in \Lambda(\Gamma_0[\mu]) : \text{dist}(z, \delta_0)\} < \epsilon.$$

*Proof.* Let $p_1, p_2, p_3 \in \Lambda(\Gamma_0[\mu_{p/q}])$ be three distinct fixed points, and $c = c(\mu)$ the center of $\delta_0$. There is a neighborhood $U \subset \mathbb{C}$ of $\mu_0$ such that the maps $p_i \colon U \to \mathbb{C}$, $i = 1, 2, 3$, and $c \colon U \to \mathbb{C}$ induced by the map $\mu \mapsto T[\mu]$, are analytic functions. Thus there is a $\eta > 0$ such that for $|\mu - \mu_0| < \eta$, $|p_i(\mu) - p_i(\mu_0)| < \epsilon/2$, $i = 1, 2, 3$, and $|c(\mu) - c(\mu_0)| < \epsilon/2$. Now for any point $z \in \Lambda(\Gamma_0[\mu])$ we have
$$|z - c(\mu_0)| = |p_1 - \mu_0|$$
$$\leq |p_1(\mu) - p_1(\mu_0)| + |p_1(\mu_0) - c(\mu_0)| + |c(\mu_0) - c(\mu)|$$
$$< |p_1(\mu_0) - c(\mu_0)| + \epsilon. \quad \square$$

**6.5. Theorem.** *On the extension $\mathscr{P}^+_{p/q}$ of each rational pleating ray there is an open neighborhood $U$ of the cusp $\mu_{p/q}$ on $\mathscr{P}^+_{p/q}$ such that the following holds: If $\mu \in \mathscr{P}^+_{p/q} \cap U$ and $|\operatorname{tr} W_{p/q}| = 2\cos(\pi/n)$ for some $n \in \mathbb{N}$, $n \geq 2$, then*
$$G[\mu] = F *_{W_{r/s}}$$
*is a Kleinian group representing a thrice punctured sphere and a sphere with a puncture and two branch points of order $n$ on its ordinary set. If $\mu \in \mathscr{P}^+_{p/q} \cap U$ is not of this form, $G[\mu]$ is not discrete.*

*Proof.* The idea of the proof is to show that we can apply Theorem 5.1 in our situation. At the cusp $\mu_{p/q} \in \partial \mathscr{M}$ the tangent combinatorial $p/q$ chain $(\Delta_i)_{i \in \mathbb{Z}}$ can be written as the collection
$$\left\{W^i_{p/q}(\Delta_j) : i \in \mathbb{Z},\ 1 \leq j \leq q-1\right\}.$$

By Proposition 6.3 we have
$$\Delta_i \cap \Delta_j = \emptyset$$
if $|i-j| > 1$. By Lemma 6.4 the finite piece $\Delta_0, \ldots, \Delta_{q_-}$ is deformed continuously in $\mu$, and no new intersections of the $\Delta_i[\mu]$ are produced if $|\mu_{p/q} - \mu|$ is small. Thus we can change $\mu$ in a small open set $U$ containing $\mu_{p/q}$ such that $(*)$ holds for all $\mu \in U$ and all $i = 0, \ldots, q$. Let $\mu_{p/q}(n)$ be a parameter on $\mathscr{P}_{p/q}^+$ such that $W_{p/q}$ is elliptic. If $\mu_{p/q}(n) \in U \cap \mathscr{P}_{p/q}^+$, then

$$\{W_{p/q}[\mu_{p/q}(n)]^i \left(\Delta_j[\mu_{p/q}(n)]\right) : i \leq i \leq n,\, 0 \leq j \leq q-1\}$$

form a finite combinatorial $-s/q$ chain. Also, after possibly making $U$ slightly smaller the disks $\mathscr{D}_A$ and $\mathscr{D}_B$ are disjoint.

In order to apply apply Theorem 5.1 on discreteness and circle chains we need that the word $W_{-s/q}[W_{p/q}[S,T], W_{r/s}[S,T]]$ is parabolic for $\mu = \mu_{p/q}(n)$. But by Proposition 6.2 this holds for any $\mu$. □

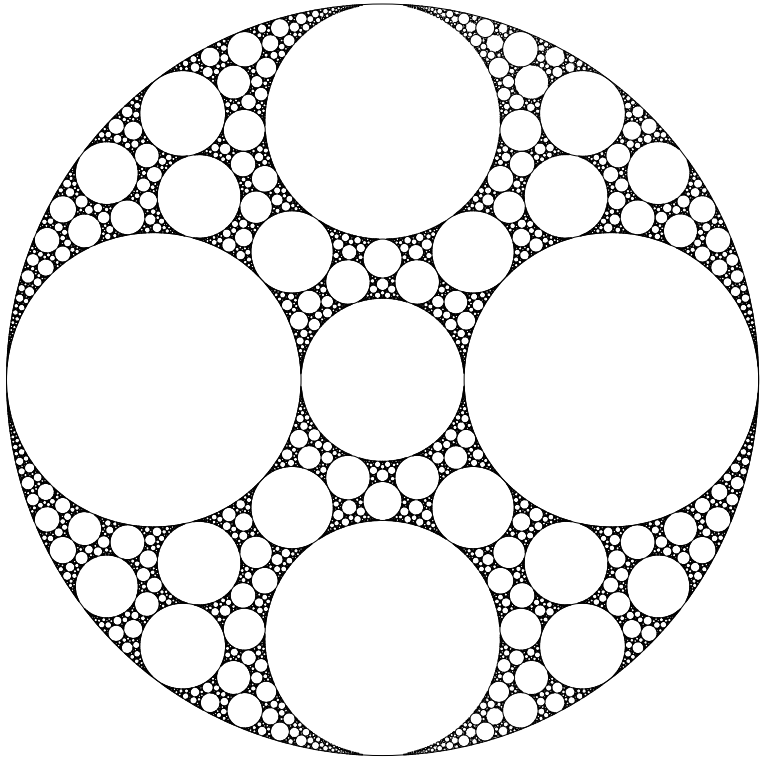

FIGURE 8. The limit set the Kleinian group on the extended ray $\mathscr{P}_{1/2}^+$ for the parameter $\mu$ such that $W_{1/2}$ is an elliptic of order 4.

We can now prove the second main result of this paper. Theorem 6.6 is a version of Theorem III of Keen, Maskit and Series [5] on the uniqueness of maximally parabolic groups applied to our situation with elliptic transformations in the groups. The proof is essentially the same:

**6.6. Theorem.** *Let $\tau_{p'/q}$ be the boundary point of $\mathscr{M}_n$ on the ray $\mathscr{P}_{p'/q}$. Then the group $G[\mu_{p/q}(n)]$ is conjugate to $G_n[\tau_{p'/q}]$ by a Möbius transformation.*

*Proof.* Both groups are geometrically finite, and they are HNN extensions of a triangle group of signature $(n, n, \infty)$. There is a natural type-preserving isomorphism $\varphi \colon G_n[\tau_{p'/q}] \to G[\mu_{p/q}(n)]$:

$$\varphi(A_n) = X, \quad \varphi(C_n[\tau]) = Z.$$

The theorem now follows from Tukia's version of Marden's Isomorphism Theorem [24] Theorem 4.2, if we can construct a conformal map

$$\psi : \Omega(G_n[\tau_{p'/q}]) \to \Omega(G[\mu_{p/q}(n)])$$

that induces $\varphi$. This can be done as follows: There is a unique Möbius transformation $\psi_1 \colon \mathbb{D} \to \delta_0$ conjugating the stabilizing triangle groups of these disks such that
$$\psi_1 \circ A_n \circ \psi_1^{-1} = X \quad \text{and} \quad \psi_1 \circ B_n \circ \psi_1^{-1} = Y.$$

We can extend this map equivariantly to the $G_n[\tau_{p'/q}]$ orbit of $\mathbb{D}$ by setting

$$\psi_1|_{w(\mathbb{D})} = \varphi(w) \circ \psi_1 \circ w^{-1}$$

for any $w \in G_n[\tau_{p'/q}(n)]$. Let $\widetilde{\delta}_0$ be the first circle of the tangent circle chain of $G_n[\tau_{p'/q}(n)]$, tangent to $\mathbb{D}$ at $\mathrm{fix}\, K_n = e^{\pi i/n}$. There is a similar map $\psi_2$ for the $G_n[\tau_{p'/q}(n)]$ orbit of $\widetilde{\delta}_0$. Together $\psi_1$ and $\psi_2$ define a conformal map of the sets of discontinuity. It is compatible with the isomorphism $\varphi$ by construction. □

## References

[1] P. Arés Gastesi, *On Teichmüller spaces of Koebe groups*, preprint (1995).
[2] A.F. Beardon, *The Geometry of Discrete Groups*, Springer-Verlag, 1983.
[3] M. Cohen, W. Metzler and A. Zimmermann, *What does a basis of $F(a,b)$ look like?*, Math. Ann. **257** (1981), 435-445.
[4] T. Jørgensen, *On discrete groups of Möbius transformations*, Amer. J. Math. **98** (1976), 739-749.
[5] L. Keen, B. Maskit and C. Series, *Geometric finiteness and uniqueness for Kleinian groups with circle packing limit sets*, J. reine angew. Math. **436** (1993), 209-219.
[6] L. Keen and C. Series, *Pleating coordinates for the Teichmüller space of a puncture torus*, Bull. Amer. Math. Soc. **26** (1992), 141-146.
[7] ———, *Pleating coordinates for the Maskit embedding of the Teichmüller space of punctured tori*, Topology **32** (1993), 719-749.
[8] ———, *The Riley slice of Schottky space*, Proc. London Math. Soc. **69** (1994), 72-90.
[9] ———, *Continuity of convex hull boundaries*, Pacific J. Math. **168** (1995), 183-206.
[10] I. Kra, *Horocyclic coordinates for Riemann surfaces and moduli spaces. I: Teichmüller and Riemann spaces of Kleinian groups*, Jour. Amer Math. Soc. **3(3)** (1990), 500-578.
[11] I. Kra and B. Maskit, *Remarks on projective structures*, Riemann Surfaces and Related Topics: Proceedings of the 1978 Stony Brook Conference, Princeton University Press, Princeton, 1980.
[12] I. Kra and B. Maskit, *The deformation space of a Kleinian group*, Amer. J. Math. **103(5)** (1980), 1065-1102.
[13] A. Marden, *The geometry of finitely generated kleinian groups*, Ann. of Math. **99** (1974), 383-462.
[14] B. Maskit, *Moduli of marked Riemann surfaces*, Bull. Amer. Math. Soc. **80** (1974), 773-777.
[15] ———, *Intersections of component subgroups of Kleinian groups*, Discontinuous groups and Riemann Surfaces: Proceedings of the 1973 Conference at the University of Maryland, Princeton University Press, Princeton, 1974.

[16] \_\_\_\_\_\_, *Kleinian Groups*, Springer-Verlag, New York, Heidelberg, Berlin, 1987.
[17] \_\_\_\_\_\_, *On Klein's combination theorem. IV*, Trans. Amer. Math. Soc. **336** (1993), 265-294.
[18] G. McShane, J. Parker and I. Redfern, *Drawing limitsets of Kleinian groups using finite state automata*, Experimental Math **3** (1994), 153-170.
[19] R. P. Osborne and H. Zieschang, *Primitives in the group on two generators*, Invent. Math **63** (1981), 17-24.
[20] J. Parkkonen, *Geometric complex analytic coordinates for deformation spaces of Koebe groups*, Ann. Acad. Sci. Fenn. Ser. I A Math. Dissertationes **102** (1995).
[21] R. Riley, *Algebra for heckoid groups*, Trans. Amer. Math. Soc. **334** (1992), 389-409.
[22] C. Series, *The geometry of Markoff numbers*, Math. Intelligencer **7(3)** (1985), 20-29.
[23] W. Thurston, *The Geometry and Topology of Three-manifolds*, lecture notes, Princeton University, 1979.
[24] P. Tukia, *On Isomorphisms of geometrically finite Möbius groups*, Publ. Math. IHES **61** (1985), 171-214.
[25] D. Wright, *The shape of the boundary of Maskit's embedding of the Teichmüller space of once-punctured tori*, preprint (1990).

DEPARTMENT OF MATHEMATICS, UNIVERSITY OF JYVÄSKYLÄ, P.O. BOX 35, FIN-40351 JYVÄSKYLÄ, FINLAND
*E-mail address*: `parkkone@math.jyu.fi`